\documentclass[12pt,a4paper]{amsart}
\usepackage{a4wide}
\usepackage[english]{babel}
\usepackage[T2A]{fontenc}
\usepackage[utf8]{inputenc} 
\usepackage{amsfonts}
\usepackage{amssymb, amsthm, amscd}
\usepackage{amsmath}
\usepackage{mathtools}
\usepackage{needspace}
\usepackage{etoolbox}
\usepackage{lipsum}
\usepackage{comment}
\usepackage{cmap}
\usepackage[pdftex]{graphicx}
\usepackage[unicode]{hyperref}
\usepackage[matrix,arrow,curve]{xy}
\usepackage[usenames,dvipsnames]{xcolor}
\usepackage{colortbl}
\usepackage{textcomp}
\usepackage{cite}
\usepackage{euscript}
\usepackage{breqn}

\newcommand{\Aut}{\,\mathrm{Aut}\,}

\newcommand{\ad}{\,\mathrm{ad}\,}

\newcommand{\GL}{\,\mathrm{GL}\,}
\newcommand{\SL}{\,\mathrm{SL}\,}
\newcommand{\SO}{\,\mathrm{SO}\,}
\newcommand{\Sp}{\,\mathrm{Sp}\,}
\newcommand{\Spin}{\,\mathrm{Spin}\,}
\newcommand{\PGL}{\,\mathrm{PGL}\,}
\newcommand{\PSL}{\,\mathrm{PSL}\,}

\newcommand{\Exp}{\,\mathrm{Exp}\,}

\newcommand{\diag}{{\rm diag}}

\newcommand{\Z}{\mathbb{Z}}
\newcommand{\N}{\mathbb{N}}

\newcommand{\F}{\mathbb{F}}

\newcommand{\ph}{\varphi}
\renewcommand{\ge}{\geqslant}
\renewcommand{\le}{\leqslant}

\renewcommand{\P}{\EuScript{P}}

\newcommand{\map}[3]{#1\colon #2\to #3}

\newcommand{\tc}{\text{,}}
\newcommand{\tp}{\text{.}}

\newcommand{\vpi}{\varpi}
\newcommand{\eps}{\varepsilon}
\newcommand{\sub}{\subseteq}

\DeclareMathOperator{\IsElUnip}{IsElUnip}
\DeclareMathOperator{\Add}{Add}
\DeclareMathOperator{\Mult}{Mult}
\DeclareMathOperator{\MatUpToCent}{MatUpToCent}
\DeclareMathOperator{\DecompConj}{DecompConj}
\DeclareMathOperator{\Matrix}{Matrix}

\theoremstyle{plain}
\newtheorem{thm}{Theorem}[section]
\newtheorem{lem}[thm]{Lemma}
\newtheorem{prop}[thm]{Proposition}
\newtheorem{cor}[thm]{Corollary}

\theoremstyle{definition}
\newtheorem{defn}[thm]{Definition}

\theoremstyle{remark}
\newtheorem{rem}[thm]{Remark}
\newtheorem{exmp}[thm]{Example}

\AtBeginEnvironment{thm}{\begin{samepage}}
	\AtEndEnvironment{thm}{\end{samepage}}
\AtBeginEnvironment{lem}{\begin{samepage}}
	\AtEndEnvironment{lem}{\end{samepage}}
\AtBeginEnvironment{defn}{\begin{samepage}}
	\AtEndEnvironment{defn}{\end{samepage}}
\AtBeginEnvironment{cor}{\begin{samepage}}
	\AtEndEnvironment{cor}{\end{samepage}}
\AtBeginEnvironment{prop}{\begin{samepage}}
	\AtEndEnvironment{prop}{\end{samepage}}

\makeatletter
\def\@settitle{\begin{center}%
		\baselineskip14\p@\relax
		\bfseries
		\@title
	\end{center}%
}

\def\@evenhead{\hfil\sc E. Bunina and P. Gvozdevsky\hfil}
\def\@oddhead{\hfil\sc Regular bi-interpretability \hfil}
\makeatother

\title{Regular bi-interpretability and finite axiomatizability of Chevalley groups}

\keywords{Chevalley groups, commutative rings, first order logic, bi-interpretability, finite axiomatizability, elementary definability}
\subjclass[2020]{20G35, 20A15, 03C60} 

\author{Elena Bunina and Pavel Gvozdevsky}
\date{}
\address{Department of Mathematics, Bar-Ilan University, 5290002 Ramat Gan, ISRAEL}
\thanks{The paper is written as part of the second author's post-doctoral fellowship at Bar-Ilan University, Department of Mathematics; and is supported by ISF grant 1994/20.}
\email{helenbunina@gmail.com, gvozdevskiy96@gmail.com}
\begin{document}

\maketitle

\begin{flushright}
    The paper is dedicated to the memory\\ of the outstanding mathematician\\ Professor Boris Plotkin
\end{flushright}

\begin{abstract}
In this paper we consider Chevalley groups over commutative rings with~$1$, constructed by irreducible root systems of rank $>1$. We always suppose that for the systems $A_2, B_\ell, C_\ell, F_4, G_2$ our rings contain $1/2$ and for the system $G_2$ also $1/3$. Under these assumptions we prove that the central quotients of Chevalley groups are regularly bi-interpretable with the corresponding rings, the class of all central quotients of Chevalley groups of a given type is elementarily definable and even finitely axiomatizable (see Definition~\ref{FAdefn}). The same holds for adjoint Chevalley groups and for bondedly generated Chevalley groups.  We also give an example of Chevalley group with infinite center, which is not bi-interpretable with the corresponding ring and is elementarily equivalent to a group that is not a Chevalley group itself. 
\end{abstract}

\section{Introduction. Historical overview}

\subsection{Elementary equivalence of groups, Maltsev-type theorems}

Elementary (first-order) classification of algebraic structures goes back to the works of Tarski and Maltsev. In general, the task is to characterize, in somewhat algebraic terms, all algebraic structures elementarily equivalent to a given one. Recall, that two algebraic structures $\mathfrak A$ and $\mathfrak B$ in a language~$\mathcal L$ are \emph{elementarily equivalent} ($\mathfrak A \equiv  \mathfrak B$) if they satisfy precisely the same
first-order sentences in~$\mathcal L$.

Classification of objects up to an equivalence relation is one of the principal problems in mathematics.
The equivalence relation on a structure is determined by the type of property one would like to be preserved,
namely in model theory one classifies modulo elementary theory; in algebra, modulo isomorphism;
in geometric group theory, modulo quasi-isometries, etc. 
To determine which groups (or other structures) are elementarily equivalent to a given one is usually a very hard problem. The first powerful  result for groups is the classification of groups elementarily equivalent
to an Abelian group. It was obtained by W.\,Szmielew~\cite{Shmelew} (two Abelian groups are elementarily equivalent if and only if they have the same special ``characteristic numbers'' like $\Exp A$, $\dim p^nA[p]$, etc.).
Nowadays the detailed descriptions of the groups elementary equivalent to a given one are obtained for wide list of groups, see for example the survey in~\cite{Kazachkov}.

When we restrict ourselves to the class of  finitely generated groups elementarily equivalent to the given one, the problem becomes especially challenging.

An outstanding result was the answer to the old problem raised by A.\,Tarski around 1945:
for non-Abelian finitely generated free groups the elementary theory doesn't distinguish these groups (see for details the series of works of Kharlampovich--Myasnikov and Z.\,Sela, e.\,g. \cite{Kharlamp-Myasikov},~\cite{Sela}). The similar situation takes place for the torsion free hyperbolic groups (see Sela,~\cite{Sela2}).

 On the other hand, a theorem of F.\,Oger~\cite{kazach091} states that two finitely generated
nilpotent groups $G$ and $H$ are elementarily equivalent if and only if $G \times \mathbb Z \cong H \times \mathbb Z$.  In the paper~\cite{M87}, A.\,Myasnikov
gave a description of groups elementarily equivalent to a given finitely generated nilpotent $K$-group,
where $K$ is a field of characteristic zero. Recently, A.\,Myasnikov and M.\,Sohrabi classified
groups elementarily equivalent to a free nilpotent $R$-group of finite rank, where $R$ is a domain of
characteristic zero, see~\cite{MS10}. M.\,Casals-Ruiz and I.\,Kazachkov in~\cite{Kazachkov}  proved that all finitely generated groups elementarily equivalent to
a solvable Baumslag--Solitar group $BS(m, 1)$ are isomorphic to it.

Certainly, any new such complete classification result is a rather rare mathematical phenomenon. 

Another way of studying elementary equivalence of algebraic models is to establish a connection between derivative models (linear groups over rings and fields, automorphism groups and endomorphism rigs of different structures, etc.) and initial models (and ``parameters'') used for the construction. In his classical paper \cite{Maltsev} A.I.\,Maltsev studied elementary equivalence of matrix groups $\mathcal G_n(F)$ where $\mathcal G_n =\GL_n, \SL_n, \PGL_n, \PSL_n$, $n\geqslant 3$, and $F$ is a field. Namely, he showed that $\mathcal G_n(F)\equiv \mathcal G_m(L)$ if and only if $n=m$ and $F\equiv L$. His proof was based on two principal results. The first one states that for any integer $k\geqslant 3$ and $\mathcal G_n$ as above there is a group sentence $\Phi_{k,\mathcal G}$ such that for any~$n$, and a field~$F$, $\Phi_{k,\mathcal G}$ holds in $\mathcal G_n(F)$ if and only if $k=n$. The second one is that $F$ and $\mathcal G_n(F)$ are mutually interpretable in each other. More precisely, $\mathcal G_n(F)$ is absolutely interpretable in~$F$ (i.\,e., no use of parameters), while $F$ is interpretable in $\mathcal G_n(F)$ uniformly with respect to some definable subset of tuples of parameters (so-called \emph{regular} interpretability). This implies that the theories $Th(F)$ and $Th(\mathcal G_n(F))$ are reducible to each other in polynomial time, hence $Th(\mathcal G_n(F))$ is decidable if and only if $Th(F)$ is decidable. Later Beidar and Mikhalev introduced another general approach to elementary equivalence of classical matrix groups~\cite{BeidarMikhalev}. Their proof was based on Keisler--Shelah theorem (\emph{two structures are elementarily equivalent if and only if their ultrapowers over non-principal ultraflters are isomorphic}, see~\cite{Keisler},~\cite{Shelah}) and the description of the abstract isomorphisms of the groups of the type $\mathcal G_n(F)$. Then these results were extended to unitary linear and Chevalley groups  (see~\cite{Bunina_intro1}, \cite{Bunina_intro3}, \cite{Bunina-local}, \cite{Bunina_intro2}).  Note that in all the results above the first-order theories include only the standard constants from the languages of groups and rings. The model theory of the group $UT_n(R)$, where $n\geqslant 3$, and $R$ is an arbitrary unitary associative ring, was studied in details by O.\,Belegradek~\cite{Bel94}. He used heavily that the ring $R$ is interpretable
(with parameters) in $UT_n(R)$. A.\,Myasnikov and M.\,Sohrabi studied model theory of groups $\SL_n(\mathcal O)$, $\GL_n(\mathcal O)$, and $T_n(\mathcal O)$ over fields and rings of algebraic integers in~\cite{MyasSohr57} and \cite{Myasnikov-Sohrabi}.
Their method exploits the mutual interpretability (and also bi-interpretability)
of the group and the ring. In a similar manner N.\,Avni, A.\,Lubotsky, and C.\,Meiri in~\cite{MyasKharl3} 
proved the first order rigidity of non-uniform higher rank arithmetic groups, and then Avni and Meiri in~\cite{AvniMeiri2} did the same  for centerless irreducible higher rank arithmetic lattices in characteristic zero. 

D.\,Segal and K.\,Tent (see~\cite{Segal-Tent}) showed that for Chevalley groups $G_{\P}(\Phi,R)$ of rank
$>1$ over an integral domain~$R$ if $G_{\P}(\Phi,R)$ has finite elementary width or is adjoint, then $G_{\P}(\Phi,R)$ and $R$
are bi-interpretable.  In~\cite{bunina2022} it was proved that over local rings Chevalley groups $G_{\P}(\Phi,R)$ of rank
$>1$ are regularly bi-interpretable with the corresponding rings.

\subsection{Chevalley groups, their automorphisms, isomorphisms and elementary equivalence.}

 In the 1950’s Chevalley, Steinberg and others introduced the concept of Chevalley groups
over commutative rings, which includes classical linear groups (special linear $\SL_n$, special orthogonal $\SO_n$, symplectic $\Sp_{2n}$,
spinor $\Spin_n$, and also projective groups connected with them).

Isomorphisms and automorphisms of Chevalley groups were also studied intensively.
The standard description of isomorphisms of Chevalley groups over fields was obtained by
R.\,Stein\-berg~\cite{Stb1} for the finite case and by J.\,Humphreys~\cite{H} for the infinite one. Many papers were devoted
to description of automorphisms of Chevalley groups over different
commutative rings (see
Borel--Tits~\cite{v22}, Carter--Chen~Yu~\cite{v24},
Chen~Yu~\cite{v25}--\cite{v29}, E.\,Abe~\cite{Abe_OSN}, A.\,Klyachko~\cite{Klyachko} and others).
Then automorphisms of arbitrary Chevalley groups over local rings were described (these rings supposed to contain~$1/2$ for the root systems $A_2, B_\ell, C_\ell$ and $F_4$ and $1/3$ for the root system~$G_2$) and then in 2012 (see~\cite{Bunina_main}) all automorphisms of adjoint Chevalley groups over arbitrary commutative rings were described. Recently (see~\cite{Bunina_all_aut})  the complete desription of automorphisms and isomorphisms of Chevalley groups over arbitrary commutative rings subject to natural invertability restrictions was obtained:

{\thm\label{Chev_Aut_all}
Let $G=G_{\P}(\Phi,R)$ be a Chevalley group  (or $G=(E_{\P}(\Phi,R)$
be its elementary subgroup) with an irreducible root system~$\Phi$ of rank $>1$,  $R$ be a commutative ring with~$1$. Suppose that for $\Phi = A_2, B_\ell, C_\ell$ or $F_4$ we have $1/2\in R$, for $\Phi=G_2$ we have $1/2,1/3 \in R$. Let us also exclude from our consideration the semi-spinor group in the case~$D_{2\ell}$. 

Then every automorphism of the group~$G$
is standard (the composition of ring, central, diagram and inner automorphisms). If $G$ is an adjoint Chevalley groups, then the inner automorphism in the composition is strictly inner and the central automorphism is trivial.
}

\medskip

Description of  automorphisms and isomorphisms of Chevalley groups over rings is a powerful instrument needed for solutions of diefferent model-theoretic problems of Chevalley groups in general (including Maltsev-type theorems, bi-interpretability, first order rigidity, elementary definability etc.). 
Talking about Maltsev-type theorems, the following result (see~\cite{Bunina_recent} and \cite{Bunina_all_aut}) almost
 finalizes  Maltsev-type theory for Chevalley groups (in some root systems we still need $1/2$ or $1/3\in R$, therefore this result is not absolutely complete):

{\thm\label{Chev_EE_all}
Let
 $G=G_{\P} (\Phi,R)$ and $G'=G_{\P'}(\Phi',R')$
\emph{(}or $E_{\P} (\Phi,R)$ and $E_{\P'}(\Phi',R'))$ be two 
\emph{(}elementary\emph{)} Chevalley groups  with indecomposable root
systems $\Phi,\Phi'$ of ranks $> 1$ over infinite commutative
rings $R$ and~$R'$  (for the cases of the roots systems $A_2$, $B_\ell$, $C_\ell$, $F_4$ with $1/2$
and for the system $G_2$ with $1/2$ and $1/3$).
 
Then the groups $G$ and $G'$ are elementarily equivalent if and only if the root systems  $\Phi$
 and $\Phi'$ are isomorphic, the weight lattices of~$\P$ and~$\P'$ are isomorphic (except the case of the root system $D_{2\ell}$ and an intermediate weight lattice), and the rings $R$  and $R'$
are elementarily equivalent.
}

\subsection{First order rigidity and QFA property, bi-interpretability of linear groups and their elementary definability.}

Recent years have been characterized by a real breakthrough in model theory of linear groups. Within quite a few years  the papers of A.\,Nies (see~\cite{Nies_main}, \cite{Nies19}, etc.), N.\,Avni, A.\,Lubotzky, C.\,Meiri (see \cite{MyasKharl3}, \cite{AvniMeiri2}, also see~\cite{MyasKharl30}), A.\,Myasnikov, O.\,Kharlampovich and M.\,Sohrabi (see~\cite{Myas-Kharl-Sohr}, \cite{Myasnikov-Sohrabi}, \cite{Myasnikov-Sohrabi2}, etc.),  
D.\,Segal and K.\,Tent (see~\cite{Segal-Tent}), B.\,Kunyavskii, E.\,Plotkin, N.\,Vavilov (see~\cite{VavKunPlot}) and others appeared (see~\cite{plotkin_last} for a survey of results). 

Since according to L\"owenheim--Skolem theorem
for any infinite structure $\mathcal A$ there exists a structure $\mathcal B$ such that $\mathcal A\equiv \mathcal B$ and they
have different cardinalities, the classification up to elementary equivalence is always strictly wider than classification up to isomorphism. But for  finitely generated  case there exists a chance that elementary equivalence could imply isomorphism.
For example if any finitely generated ring $A$ is elementarily equivalent to~$\mathbb Z$, then it is isomorphic to~$\mathbb Z$. This phenomenon is not true for infinitely generated rings.

A question dominating research in this area has been if and when elementary
equivalence between finitely generated groups (rings) implies isomorphism. Recently,
Avni--Lubotzky--Meiri~\cite{MyasKharl3}  found the term \emph{first-order rigidity}.
Namely, a finitely generated
group (ring or other structure) $\mathcal A$ is \emph{first-order rigid} if any other finitely generated group (ring or other structure)
elementarily equivalent to~$\mathcal A$ is isomorphic to~$\mathcal A$.

In many cases a more strong QFA property holds. It says that 
 a single first-order axiom is sufficient to distinguish a group among
all the finitely generated groups. The corresponding concept was introduced in \cite{Nies19} for
groups only, and then in \cite{Nies_main} for arbitrary structures.

For a fixed finite signature an infinite finitely generated structure $\mathcal A$ is \emph{quasi finitely
axiomatizable} (\emph{QFA}) if there is a first order sentence $\varphi$ such that

(1) $\mathcal A\vDash \varphi$;

(2)  if $H$ is a finitely generated structure in the same signature such that $H \vDash \varphi$, then
$H\cong \mathcal A$.

The most important example of QFA rings is the ring of integers $(\mathbb Z, +, \cdot)$, see~\cite{Sabbagh}.

\smallskip

New horizons in this area were opened towards application of the bi-interpretability property and finding  new classes of QFA rings and groups.

{\defn\label{bi-def}
Suppose structures $\mathcal A, \mathcal B$ in finite signatures are given, as well as
interpretations of $\mathcal A$ in~$\mathcal B$, and vice versa. Then an isomorphic copy~$\widetilde{\mathcal A}$ of~$\mathcal A$ can be defined in~$\mathcal A$, by ``decoding''~$\mathcal A$ from the copy of~$\mathcal B$ defined in~$\mathcal A$.
 Similarly, an isomorphic copy $\widetilde{\mathcal B}$ of~$\mathcal B$ can be defined in~$\mathcal B$. An isomorphism $\Phi: {\mathcal A}\cong \widetilde{\mathcal A}$ can be viewed as a relation on~$\mathcal A$, and similarly for an isomorphism $\mathcal B \cong \widetilde{\mathcal B}$. 
We say that $\mathcal A$ and $\mathcal B$ are (\emph{weakly}) \emph{bi-interpretable} (with parameters) if there
are such isomorphisms that are first-order definable (see~\cite{Nies7}, Ch.\,5).
}

\medskip

The case of  special interest is when  one can prove that a given structure is bi-interpretable with the ring $\mathbb Z$ of integers.
A.\,Khelif (see~\cite{Nies12}, see also~\cite{Myas-Kharl-Sohr}) realized that one can use bi-interpretability of a finitely generated structure~$\mathcal A$ with~$(\mathbb Z, + \times)$ as a general method
to prove that $\mathcal A$ is QFA. 
Somewhat later, P.\,Dittmann and  F.\,Pop (see~\cite{DittmannPop}) used this method to show that each finitely generated  field is QFA. 

Anyway, suppose the structure $\mathcal A$ in a finite signature is bi-interpretable with the ring~$\mathbb Z$. Then $\mathcal A$ is QFA.

\smallskip

Along with  Maltsev-type theorems, the problems of the \emph{elementary definability} of classes of groups are quite important. These problems are formulated as follows: 
\emph{Suppose that we have a class $\mathcal G$ of groups \emph{(}ring, etc.\emph{)} and an arbitrary group \emph{(}rings, etc.\emph{)} $H$ which is elementarily equivalent to some $G\in \mathcal G$. Is it true that $H\in \mathcal G$?}

For example  the class of all linear groups over fields (commutative rings) is elementarily definable (A.I.\,Maltsev,~\cite{Mal}),
but on the other hand for general linear groups even over fields there is no elementary definability: there exists a group $H\equiv \GL_n(F)$ which is not of type $\GL$ itself.

In 1984 B.I.\,Zilber (see~\cite{Zilber1984}) proved that for algebraic groups over algebraically closed fields  the problem of elementary definability has the positive answer:
\emph{the class of all simple algebraic groups over algebraically closed fields  is elementarily definable.}

To extend this Zilber's result for the case of rings  it is crucial to introduce a notion of  regular bi-interpretability.
If $\mathcal A$ is interpreted in~$\mathcal B$ uniformly with respect
to a $\emptyset$-definable subset $D\subseteq \mathcal B^k$ then we say that $\mathcal A$ is \emph{regularly interpretable} in~$\mathcal B$ and write in this case $\mathcal A \cong \Gamma(\mathcal B, \varphi)$, provided $D$ is defined by~$\varphi$ in~$\mathcal B$ (see~\cite{DaniyarMyasn}).

{\defn\label{reg-bi-def}[see \cite{DaniyarMyasn}]
Algebraic structures $\mathcal A$ and $\mathcal B$ are called (\emph{strongly}) {\em regularly bi-interpretable}, if  
\begin{enumerate}
\item there exist a regular interpretation $(\Gamma,\varphi)$ of $\mathcal A$ in $\mathcal B$ and a regular interpretation $(\Delta,\psi)$ of $\mathcal B$ in $\mathcal A$;
\item there exists formula $\theta_{\mathcal A}(\bar u, x, \bar r)$ in $L(\mathcal A)$, where $|\bar u|={\dim\Gamma\cdot\dim\Delta}$, $|\bar r|={\dim_{par}\Gamma\circ\Delta}$, such that for any tuple $\bar r_0\in \varphi_\Delta\wedge \psi(\mathcal A)$ the formula $\theta_{\mathcal A}(\bar u, x, \bar r_0)$ defines some coordinate map $U_{\Gamma\circ\Delta}(\mathcal A,\bar r_0)\to A$;
\item there exists formula $\theta_{\mathcal B}(\bar u, x, \bar t)$ in $L(\mathcal B)$, where $|\bar u|={\dim\Gamma\cdot\dim\Delta}$, $|\bar t|={\dim_{par}\Delta\circ\Gamma}$, such that for any tuple $\bar t_0\in \psi_\Gamma\wedge \varphi(\mathcal B)$ the formula $\theta_{\mathcal B}(\bar u, x, \bar t_0)$ defines some coordinate map $U_{\Delta\circ\Gamma}(\mathcal B,\bar t_0)\to B$. 
\item for any 
pair of parameters $(\bar p, \bar q)$, $\bar p\in \varphi(\mathcal B)$, $\bar q\in \psi(\mathcal A)$, there exists a pair of coordinate maps $(\mu_\Gamma,\mu_\Delta)$ for interpretations $(\Gamma,\bar p)$ and $(\Delta,\bar q)$, such that for any $\bar r_0=(\bar{\bar p},\bar q)$, $\bar{\bar p}\in \mu^{-1}_\Delta(\bar p)$, and $\bar t_0=(\bar{\bar q},\bar p)$, $\bar{\bar q}\in \mu^{-1}_\Gamma(\bar q)$, the coordinate maps ${\mu_\Gamma\circ\mu_\Delta\colon} U_{\Gamma\circ\Delta}(\mathcal A,\bar r_0)\to A$ and ${\mu_\Delta\circ\mu_\Gamma\colon} U_{\Delta\circ\Gamma}(\mathcal B,\bar t_0)\to B$ are defined in $\mathcal A$ and $\mathcal B$ correspondingly by the formulas $\theta_{\mathcal A}(\bar u, x, \bar r_0)$ and $\theta_{\mathcal B}(\bar u, x, \bar t_0)$.
\end{enumerate}
}

\medskip

Throughout this paper, all bi-interpretations under consideration will be both regular and strong.

%
%

\smallskip

The existence of regular bi-interpretation plays the decisive role for establishing elementary  definability of given classes of groups.
\emph{If linear/Chevalley groups (or another derivative structures) of any concrete type over some classes of fields/rings are regularly bi-interpretable with the corresponding rings, then this class of groups (structures) is elementarily definable} (see \cite{bunina2022}).

 A.\,Myasnikov and M.\,Sohrabi proved  that for all linear groups  considered in the paper~\cite{Myasnikov-Sohrabi}   the bi-interpretation with the corresponding rings is indeed a regular bi-interpretation,
hence  the elementary definability problem has positive solution. One can show that this is  also the case for  the results of Segal-Tent (see~\cite{Segal-Tent}) mentioned above.

\bigskip

\section{Main results}\leavevmode

All mentioned  recent results illustrate that along with bi-interpretability with parameters it is of unprecedented importance to state {\bf regular} bi-interpretability of linear groups and the corresponding rings.

In~\cite{bunina2022}  the following theorem was proved:

{\thm\label{th-reg-bi-interpret}
If  $G(R)=G_\P (\Phi,R)$  is a Chevalley group $($or $G(R)=E_{\P}(\Phi, R)$ is its elementary subgroup$)$ of rank $> 1$, $R$ is a local ring \emph{(}with $\frac{1}{2}$ for the root systems ${B}_\ell, {C}_\ell, {F}_4, {G}_2$ and with  $\frac{1}{3}$ for  ${G}_{2})$, then the group $G(R)$ is regularly bi-interpretable with~$R$. Therefore the class of Chevalley groups over local rings \emph{(}with mentioned restrictions\emph{)} is elementarily definable.
}

\medskip

This paper is devoted to obtaining the most general answer to the following fundamental question: 

\medskip

{\bf Problem $(*)$.} \emph{For which rings and under which conditions the class of all Chevalley groups is elementarily definable?}

\medskip

Also, we will consider  other model-theoretic questions of Chevalley groups in maximal generality.

\medskip

As it was observed in Section~1, Problem~$(*)$ has a positive solution for the Chevalley groups  $G_\P(\Phi,R)$, when $R$ is an algebraically closed field~\cite{Zilber1984}, when $R$ is a Dedekind ring of the arithmetic type (can be deduced from~\cite{Segal-Tent}), and when $R$ is a local ring~\cite{bunina2022}.  In the later two  cases some mild invertability conditions on~$R$ are  assumed.   

 The best situation would be, if for all commutative rings with~$1$ the corresponding Chevalley groups are regularly bi-interpretable with the rings, which implies the fact that Chevalley groups constitute an elementary class of groups.  Unfortunately this is not the case. We will see in Section~6 that if a Chevalley group has an infinite center, then in certain cases it is possible to prove that this group is not bi-interpretable with the corresponding ring and not elementarily definable in the class of all groups. To the contrary, if the center is trivial (it holds for all adjoint Chevalley groups over all rings), then regular bi-interpretability and elementary definability take place for ``almost all'' commutative rings, i.\,e., for all commutative rings with~$\frac{1}{2}$ for the root systems ${B}_\ell, {C}_\ell, {F}_4, {G}_2$ and with  $\frac{1}{3}$ for~${G}_{2}$. 

For  the case of non-trivial finite center the answer for  the
Chevalley groups over arbitrary rings is not known. 
However, we prove that in case of boundedly generated Chevalley groups the answer is positive (Theorem~\ref{RegularBoundGen}). The same is true for the arbitrary central quotients of Chevalley groups (Theorem~\ref{RegularTrivCent}).

Moreover, in all positive cases we show that the classes of the corresponding Chevalley groups (of a given type and with a given weight lattices) are  finitely axiomatizable in the sense of the following definition.  

{\defn\label{FAdefn} Let $S$ be a class of structures in some fixed signature. We say that $S$ is finitely axiomatizable if there is a proposition $\varphi$ in that signature such that for any structure $\mathcal{A}$ in that signature we have $\mathcal A\vDash \varphi$ if and only if $\mathcal A\in S$.}

\smallskip

To formulate short informal and detailed formal versions of two main theorems and further conjectures we need some auxiliary notations/definitions.

{\defn
For any group~$G$ we will denote by $PG$ the quotient group by its center:
$$
PG: = G/Z(G).
$$
In particular, we will use the notation $PG_\P (\Phi,R)$ and will call this group a \emph{central quotient of the Chevalley group $G_{\P}(\Phi,R)$.}
}

\smallskip

Over an algebraically closed ring $PG_\P (\Phi,R)$ coincides with the adjoint Chevalley group $G_{\ad}(\Phi,R)$, but in the general case it is not always true. To be precise the functor $G_{\ad}(\Phi,-)$ is a sheafification of  $PG_\P (\Phi,-)$ in the fppf-topology. Thus $PG_\P (\Phi,R)$ is a normal subgroup of $G_{\ad}(\Phi,R)$ and the quotient group can be embedded into the cohomology group $H^1_{\mathrm{fppf}}(R,\mu)$, where $\mu=\mathrm{Ker}(G_{\P}(\Phi,-)\to G_\P(\Phi,-))$, see, for example, \cite[Chapter III, Lemma 2.6.1]{Knus}.

{\defn
A commutative ring $R$ with~$1$ is called \emph{good}, if

--- when we deal with the root system $\Phi=A_2, B_\ell, C_\ell, F_4$, then $1/2\in R$;

--- when we deal with the root system $\Phi=G_2$, then $1/2$ and $1/3\in R$.

For the cases $\Phi=A_\ell, D_\ell, E_\ell$, $l\geqslant 3$, there are no restrictions on~$R$.
}

\smallskip

{\defn Let $G$ be a group, $X\sub G$ be a subset, and $N$ be a non-negative integer. We say that $G$ is $N$-{\it boundedly generated} by the set $X$, if every element of $G$ is the product of at most $N$ elements from $X$.}

\smallskip

{\defn\label{Xg} For a group $G$, elements $g_1$,$\ldots$,$g_l\in G$ and a natural number~$m$ we set
$$
X^{(m)}_{g_1,\ldots,g_l}=\{[a,b]\mid a,b\in G\}\cup\{a^m\mid a\in G\}\cup\{g_1,\ldots,g_l\}.
$$} 

\medskip

Now we are in a position to give a brief account of our results. 

\smallskip

Let us always assume that

--- a root system $\Phi$ is irreducible and it's rank is $>1$;

--- a ring $R$ is a good commutative ring with~$1$.

\smallskip

Also let us call  results about bi-interpretability and elementary definability \emph{positive}, if:

--- the Chevalley group $G_\P (\Phi,R)$ (or $PG_\P (\Phi,R)$) is regularly bi-interpretable with~$R$;

--- any group elementarily equivalent to $G_\P (\Phi,R)$ (or $PG_\P (\Phi,R)$) has the form $G_\P(\Phi, S)$ (resp. $PG_\P (\Phi,S)$), where the ring $S$ is elementarily equivalent to $R$;

--- the class of all Chevalley groups $G_\P (\Phi, -)$ (or $PG_\P (\Phi,-)$) over all good rings is finitely axiomatizable.

\smallskip

So, shortly speaking, in this paper we obtain {\bf positive results} for the following groups and rings:

\smallskip

{\bf 1.} All central quotients $PG_\P (\Phi,R)$ of Chevalley groups over good rings;

\smallskip

{\bf 2.} In particular, all adjoint Chevalley groups $G_\P (\Phi,R)$
over good rings;

\smallskip

{\bf 3.} All Chevalley groups $G_\P (\Phi,R)$ over good rings that are $N$-boundedly  generated by some set $X^{(m)}_{g_1,\ldots,g_l}$ for some natural~$N$ and some special~$m$.

\smallskip

Also we show the {\bf negative} result (by negative we mean that all {\bf 1}--{\bf 3} do not hold) for some $G_\P (\Phi,R)$ in the case of infinite center. 

Our negative conjectures for today are the following:

\smallskip

{\bf 1.} If the center of $G_\P (\Phi,R)$ is infinite, usually the results are negative.

\smallskip

{\bf 2.} Even if the center of $G=G_\P (\Phi,R)$ is finite and $G$ is not $N$-boundedly generated for any  $X^{(m)}_{g_1,\ldots,g_l}$ and any $N\in \mathbb N$, then there exist negative examples.

\medskip

Below are the precise formulations of two main theorems of the  paper. The first one is about central quotients of Chevalley groups: 

{\thm\label{RegularTrivCent} 
Let $\Phi$ be an irreducible root system of rank $>1$, $\P$ be any corresponding weight lattice, $G_\P(\Phi,R)$ the corresponding Chevalley group. Then there exists a $|\Phi|$-ary logical formula $ \ph(g_{\alpha}\mid \alpha\in\Phi)$ in the first order language of groups such that the following holds:
\begin{enumerate}
    \item If $R$ is a good commutative ring, then for any given elements $\{ g_{\alpha}\in PG_{\P}(\Phi,R)\mid \alpha\in\Phi\}$, we have  
$$
PG_{\P}(\Phi,R)\models  \ph(g_{\alpha}\mid \alpha\in\Phi)
$$
if and only if there exists an automorphism $\overline f\in \Aut (PG_{\P}(\Phi,R))$ such that $\overline f(g_{\alpha})=\overline{x_{\alpha}}(1)$ for all $\alpha\in\Phi$.

    \item Under the assumptions of item~\emph{1}, there is a regular bi-interpretation of $R$ and $PG_{\P}(\Phi,R)$ such that the interpretation of $PG_{\P}(\Phi,R)$ in $R$ is absolute and the parameters for interpretation of $R$ in $PG_{\P}(\Phi,R)$ are defined by~$\ph$.

    \item If $G$ is a group and $g_{\alpha}\in G$, $\alpha\in\Phi$ are such that $G\models \ph(g_{\alpha}\mid \alpha\in\Phi)$, then there exists a good ring $R$ and an isomorphism $\map{f}{G}{PG_{\P}(\Phi,R)}$ such that $f(g_{\alpha})=\overline{x_{\alpha}}(1)$ for all $\alpha\in\Phi$.

    \item If $G$ is a group, $G\equiv G_\P (\Phi,R)$, where $R$ is good, then their centers are elementarily equivalent and there exists a ring $S\equiv R$ such that
     $$
      G/Z(G) \cong PG_\P (\Phi,S).
     $$
\end{enumerate}
}

\medskip
 
The second main theorem is about Chevalley groups themselves with the property of $N$-bounded generation. Morally it says that Chevalley groups that are boundedly generated by some set $X^{(m)}_{g_1,\ldots,g_l}$ (see Definition~\ref{Xg}) have the same properties that was established for groups $PG_{\P}(\Phi,R)$ in Theorem~\ref{RegularTrivCent}, However, in this case there must be some additional parameters in the formula $\ph$, the meaning of which will become clear after formulating item~1.

Recall that Chevalley group can be defined as a subgroup of $\GL(n,R)$ that consists of matrices that satisfy certain polynomial equations with integral coefficients. In other words, there is a representation $\map{\vpi}{G_{\P}(\Phi,-)}{\GL(n,-)}$ of the Chevalley-Demazure scheme, that is a closed embedding (see~\cite{ChevalleySemiSimp}).

{\thm\label{RegularBoundGen} 
Let $\Phi$ be an irreducible root system of rank $>1$, $\P$ be any corresponding weight lattice. Let $m$ be such number that the center of the group scheme $G_{\P}(\Phi,-)$ is the group $\mu_m$ of the $m$-th roots of unity. Let $N$, $l$ and $n$ be some  non-negative integers. Let $\alpha_0\in\Phi$ be some fixed root. Let $\map{\vpi}{G_{\P}(\Phi,-)}{\GL(n,-)}$ be some faithful representation, i.e a closed embedding of algebraic group schemes. Then there exists a $|\Phi|+l(n^2+1)$-ary logical formula $\ph(g_{\alpha}\mid \alpha\in\Phi;\; g_k\mid 1\le k\le l;\; g_{ij}^{(k)}\mid 1\le k\le l,\; 1\le i,j\le n)$ in the first order language of groups such that the following holds:
\begin{enumerate}
    \item If $R$ is a good commutative ring, then for any given elements $\{ g_{\alpha}\mid \alpha\in\Phi;\; g_k\mid 1\le k\le l;\; g_{ij}^{(k)}\mid 1\le k\le l,\; 1\le i,j\le n\}$, we have  
$$
G_{\P}(\Phi,R)\models  \ph(g_{\alpha}\mid \alpha\in\Phi;\; g_k\mid 1\le k\le l;\; g_{ij}^{(k)}\mid 1\le k\le l,\; 1\le i,j\le n)
$$
if and only if the group $G_{\P}(\Phi,R)$ is $N$-boundedly generated by the set $X^{(m)}_{g_1,\ldots,g_l}$, and there exists an automorphism $f\in \Aut (G_{\P}(\Phi,R))$ such that

$\bullet$ $f(g_{\alpha})= x_{\alpha}(1)$ for all $\alpha\in\Phi$; 

$\bullet$ $f(g_{ij}^{(k)})=x_{\alpha_0}(\xi_{ij}^{(k)})$ for all $1\le k\le l;\; 1\le i,j\le n$ for some $\xi_{ij}^{(k)}\in R$;

$\bullet$ $\vpi(f(g_k))$ is the matrix $(\xi_{ij}^{(k)})$ for all $1\le k\le l$.

    \item  If $R$ is a good commutative ring, and there is a $l$-tuple $g_1,\ldots,g_l\in G_{\P}(\Phi,R)$ such that the group $G_{\P}(\Phi,R))$ is $N$-boundedly generated by the set $X^{(m)}_{g_1,\ldots,g_l}$, then there is a regular bi-interpretation of $R$ and $G_{\P}(\Phi,R)$ such that the interpretation of $G_{\P}(\Phi,R)$ in $R$ is absolute and the parameters for interpretation of $R$ in $G_{\P}(\Phi,R)$ are defined by~$\ph$.

    \item If $G$ is a group and $g_{\alpha}, g_k, g_{ij}^{(k)}\in G$, $\alpha\in\Phi,\; 1\le k\le l$, $1\le i,j\le n$, are such that $G\models \ph(g_{\alpha}\mid \alpha\in\Phi;\; g_k\mid 1\le k\le l;\; g_{ij}^{(k)}\mid 1\le k\le l,\; 1\le i,j\le n)$, then there exists a good ring $R$ and an isomorphism $\map{f}{G}{G_{\P}(\Phi,R)}$ such that 
    
    $\bullet$ $f(g_{\alpha})= x_{\alpha}(1)$ for all $\alpha\in\Phi$; 

    $\bullet$ $f(g_{ij}^{(k)})=x_{\alpha_0}(\xi_{ij}^{(k)})$ for all $1\le k\le l$ for some $\xi_{ij}^{(k)}\in R$;

    $\bullet$ $\vpi(f(g_k))$ is the matrix $(\xi_{ij}^{(k)})$ for all $1\le k\le l$.

    \item Let $G$ be a group, and $G\equiv G_\P (\Phi,R)$, where $R$ is good. Assume that for some $g_1,\ldots,g_l\in G_{\P}(\Phi,R))$ the group $G_{\P}(\Phi,R))$ is $N$-boundedly generated by the set $X^{(m)}_{g_1,\ldots,g_l}$, then there exists a ring $S\equiv R$ such that
     $$
      G \cong G_\P (\Phi,S).
     $$
\end{enumerate}
}
{\rem The statement of Theorem~\ref{RegularBoundGen} makes it natural to ask for what rings the group $G_\P(\Phi,R)$ is boundedly generated by the set $X^{(m)}_{g_1,\ldots,g_l}$ as above. Let us list what we know.

$\bullet$ If $R$ is a Dedekind domain of arithmetic type then $G_\P(\Phi,R)$ s boundedly generated by elementary unipotents, and hence by commutators, see for example~\cite{CarterKeller},\cite{Morris},\cite{VavKunPlot}.

$\bullet$ Assume that $G_\P(\Phi,R)=G_{\mathrm{sc}}(\Phi,R)$ is simply connected, and the ring $R$ is finitely generated, then by~\cite{HazVav} the group $K_1(\Phi,R)=G_{\mathrm{sc}}(\Phi,R)/E_{\mathrm{sc}}(\Phi,R)$ is nilpotent. Further assume that $K_1(\Phi,R)$ is finitely generated (the question whether this groups is always finitely generated for a finitely generated ring $R$ is an old open problem). Then $K_1(\Phi,R)$ is polycyclic, see~\cite[Chapter 4A, Exercise 1.2]{Polycyclic}. It then follows that $K_1(\Phi,R)$ is boundedly generated by $m$-th powers and the set of its generators,  which is finite. Indeed, the group $\Z$ is boundedly generated by $m$-th powers and a generator, because $\Z/m\Z$ is finite; now for a polycyclic group it follows by induction. Thus we reduce the question to an elementary group. Similarly, we can do this for a non simply connected group $G_\P(\Phi,R)$, if we assume additionally that its quotient by the image of $G_{\mathrm{sc}}(\Phi,R)$ is finitely generated; this quotient is a subgroup in the cohomology group $H^1_{\mathrm{fppf}}(R,\mu)$, where $\mu=\mathrm{Ker}(G_{\mathrm{sc}}(\Phi,-)\to G_\P(\Phi,-))$, see, for example, \cite[Chapter III, Lemma 2.6.1]{Knus}.

$\bullet$ In papers~\cite{VasBounded},\cite{GvozBoundedOrt},\cite{GvozBoundedE} it is proved that an element from the Chevalley group of type $\Phi$, where $\mathrm{rk}\, \Phi\ge 2$ and $\Phi\ne F_4,G_2,E_8$, over a polynomial ring with coefficients in a small-dimensional ring can be reduced to an element of certain proper subsystem subgroup by a bounded number of elementary root elements (the bound depends on the number of variables). By induction, for the ring $\Z[x_1,\ldots,x_n]$ we can reduce such an element to a subsystem subgroup of type~$A_1$. Since every finitely generated ring is a quotient of $\Z[x_1,\ldots,x_n]$, and the map on corresponding elementary subgroups is surjective, we conclude that for a finitely generated ring any element from the elementary subgroup can be reduced to a subsystem subgroup of type~$A_1$. Note that the systems $\Phi\ne F_4,G_2,E_8$ are inessential for us here, because for them the all Chevalley groups are adjoint; hence we can apply Theorem~\ref{RegularTrivCent} to them, so we don't need bounded generation.

$\bullet$ The question whether the group $G(\Phi,\Z[x_1,\ldots,x_n])$ is boundedly generated by elementary root elements is an open problem. The educated guesses of the specialists are that the answer is negative. Now it is natural to formulate a weaker version of that problem: {\it is the group $G(\Phi,\Z[x_1,\ldots,x_n])$ boundedly generated by commutators, $m$-th powers, and a finite set, where $m$ is as above?}}

\section{Basic notations and useful facts about Chevalley groups}

We refer for the facts related to root systems and semisimple Lie algebras to \cite{Bourbaki4-6} and~\cite{Hamfris}.
More detailed information about elementary Chevalley groups is
contained in the books~\cite{Steinberg}, and about the Chevalley
groups  (also over rings) in~\cite{VavChevalley}, \cite{VavPlotkin1}, \cite{ChevalleySemiSimp} (see also
later references in these papers).

We always fix some arbitrary irreducible root system $\Phi$ of the
rank~$\ell> 1$. 

We consider an arbitrary Chevalley group $G_{\P}(\Phi,R)$,
constructed by the root system $\Phi$, a ring~$R$ and a lattice $\P$ of a
representation~$\pi$ of the corresponding Lie algebra. It is known,
that Chevalley group is defined by these three parameters.  If we consider an elementary
Chevalley group, we denote it by $E_{\P}(\Phi,R)$.

Relations between Chevalley groups and the corresponding elementary
subgroups is an important problem in the theory of Chevalley
groups over rings. For elementary Chevalley groups there exists a
convenient system of generators $x_\alpha (\xi)$, $\alpha\in \Phi$,
$\xi\in R$, and all relations between these generators are well-known.
For general Chevalley groups it is not always true.

If $R$ is an algebraically closed field, then
$$
G_{\P} (\Phi,R)=E_{\P} (\Phi,R)
$$
for any weight lattice~$\P$. This equality is not true even for the
case of fields, which are not algebraically closed. 

If $\Phi$ is an irreducible root system of  rank $\ell>
1$, then $E_{\P}(\Phi,R)$ is always normal in $G_{\P}(\Phi,R)$ (see~\cite{Taddei}, \cite{HazVav}) and even is characteristic in it. 

We will use the following result about normal subgroups of  elementary Chevalley groups (see\cite{AbeSuzuki}).

Note that for every ideal $I$ of~$R$ the natural mapping $R \to R/I$ induces a homomorphism
$$
\lambda_I : G_\P(\Phi ,R) \to  G_\P(\Phi,R/I).
$$
If $I$ is a proper ideal of~$R$, then the kernel of~$\lambda_I$ is a non-central normal subgroup of $G_\P(\Phi,R)$.

{\defn 
We denote the inverse image of the center of $G_\P(R/I)$ under $\lambda_I$ by $Z_\P(\Phi,R, I)$.

By $E_\P(\Phi,R, I)$ we denote the minimal normal subgroup of $E_\P(\Phi,R)$ which contains all
$x_\alpha(t)$, $\alpha\in \Phi$, $t\in I$. 
}

{\thm \label{NormalStructure} 
Let the rank of an irreducible root system $\Phi$ is greater than~$1$. If a
subgroup $H$ of $E_\P(\Phi,R)$ is normal in $E_\P(\Phi,R)$, then
$$
E_\P(\Phi,R, I) \leqslant H \leqslant Z_\P(\Phi,R, I) \cap E_\P(\Phi,R)
$$
for some uniquely defined ideal $I$ of the ring~$R$.
}

For adjoint Chevalley groups (which always have trivial center) the group $Z_\P (\Phi,R,I)$ coincides with the group $G_\P (\Phi,R,I)$ consisting of all matrices $A$ from $G_\P (\Phi,R)$ such that $A-E$ have coefficients from~$I$.

To prove our main theorem, we will need to study isomorphisms between central quotients of Chevalley groups. This result is of interest by itself. 

{\thm \label{isomorphic}
If $\Phi$, $\Psi$ are irreducible root systems of ranks $>1$, $\P_1$, $\P_2$ are any corresponding lattices, $R$ and $S$ are good commutative rings with~$1$, then if the groups
$$
PG_{\P_1} (\Phi,R)=G_{\P_1} (\Phi,R)/Z(G_{\P_1} (\Phi,R))\text{ and }
PG_{\P_2} (\Psi,S)=G_{\P_2} (\Psi,S)/Z(G_{\P_2} (\Psi,S))
$$
are isomorphic, then the rings $R$ and $S$ are isomorphic and the systems $\Phi$ and $\Psi$ coincide.
}

\begin{proof}
Since the (elementary Chevalley) subgroups $E_{\P_1} (\Phi,R)$ and $E_{\P_2} (\Psi,S)$ are normal in the corresponding Chevalley groups $G_{\P_1} (\Phi,R)$ and $G_{\P_2} (\Psi,S)$, then their quotients by their centers are normal subgroups in $PG_{\P_1} (\Phi,R)$ and $PG_{\P_2} (\Psi,S)$. Since for any elementary Chevalley group $E_\P (\Phi,R)$ the quotient group by its center is elementary adjoint Chevalley group $E_{ad}(\Phi,R)$, then we can assume that
 $$
 E_{ad}(\Phi,R)\vartriangleleft PG_{\P_1} (\Phi,R)\text{ and }
 E_{ad}(\Psi,S)\vartriangleleft PG_{\P_2} (\Psi,S).
 $$

 Now we want to prove that under our isomorphism $f$ the subgroup $E_{ad}(\Phi,R)$ is mapped onto the subgroup $E_{ad}(\Psi,S)$. Clearly it is enough for our result.

 Let us prove that the subgroup $E_{ad}(\Phi,R)$ 
 can be characterized as the smallest by inclusion among all the subgroups $K \leqslant G= PG_{\P_1} (\Phi,R)$ that satisfy the following properties:

1) $K$ is normal and is generated as a normal subgroup by a single element;

2) $K = [K, K]$;

3) the centralizer of $K$ in $G$ is  trivial.

 For the elementary subgroup $E=E_{ad}(\Phi,R)$ all these three properties over good rings are well known. Now let $K \leqslant  G$ be a subgroup that satisfies (1)--(3); we must prove that $E=E_{ad}(\Phi,R)\leqslant K$.
 
Since by (1) $K$ is normal in~$G$, it follows by description of the normal subgroups that there exists a unique 
ideal $J$ of $R$ such that
$$
E_{\ad}(\Phi,R,J)\subseteq K \subseteq PG_{\P_1}(\Phi,R, J).
$$

We claim that $J = R$. By (1) $K$ is generated as a normal subgroup by a single element~$g_0$.

Now, clearly $J$ is the smallest by inclusion ideal such that $g_0 \in PG_{\P_1}(\Phi,R, J)$.

hence,
$J$ is generated by all the entries of any representor of the element $g_0-I$ in the matrix ring of the representation of the initial Chevalley group.
 Therefore, the ideal $J$ is finitely generated. Since by (2) the group $K$ is perfect, the
uniqueness of the ideal $J$ implies that $J = JJ$. By Nakayama’s Lemma, there exists $s \in R$
such that $s \equiv  1 (\mod J)$ and $sJ = 0$. Thus, $E_\ad(\Phi, sR)$ is contained in the centralizer of $K$ in~$G$;
hence, (3) implies that $s = 0$; hence, we have $J = R$; hence, we have $E_\ad (\Phi, R)\subseteq K$.

 \end{proof}

\section{Proof of Theorem~\ref{RegularTrivCent}}

First of all we will give the main scheme of the proof in simple words, and after that the detailed proof will follow.

We need in general to describe all central quotients of Chevalley groups of some given type $PG_\P (\Phi,-)$ over all good rings by one sentence.

Let us suppose that the set of parameters $M=\{ \overline{x_\alpha} (1) \mid \alpha \in \Phi\}$ is fixed. Then (see Lemma~\ref{ElementaryRootSubgDefine}) for any $\alpha\in \Phi$ the subgroup $X_\alpha:=\{ \overline{x_\alpha} (t)\mid t\in R\}$ is definable with the parameters from~$M$ (and the defining formula is the same for all rings). According to Lemma~\ref{AddAndMult} on each this set we can (by the same formulas for all rings~$R$) define operations of addition and multiplication so that $X_\alpha$ with these operations becomes isomorphic to the corresponding ring~$R$.

We deal with some algebraic group scheme and therefore for our type $G_\P (\Phi,-)$ we can fix the dimension $n$ of the corresponding $\GL_n(-)$ and the set of the corresponding algebraic equations. 

Therefore the elements of $G_\P (\Phi, R)$ can be defined as $n\times n$-matrices ($n^2$-tuples) of elements of~$R$ satisfying  given equations. 

These arguments give us mutual interpretability of $G_\P (\Phi,R)$ or $PG_\P (\Phi,R)$ and~$R$ (given by same code for all rings). 
To show bi-interpretability between $G_\P (\Phi,R)$ and~$R$ we need two more things: 

(1) Present a formula stating an isomorphism between~$R$ and the matrices representing~$R$ in $PG_\P (\Phi,R)$, i.\,e., the matrices $\overline{x_\alpha} (t)$ for some fixed~$\alpha$. It is easy, since we know the form of $x_\alpha (t)$ in $\GL_n(R)$. The similar formula is of course possible for the Chevalley group $G_\P (\Phi,R)$ itself, not only for its central quotient.

(2) Present a formula stating an isomorphism between the group $PG_\P (\Phi,R)$ (as an abstract group) and the matrices from $\GL_n(R)$ satisfying the corresponding relations (see Lemma~\ref{MatrixUpToCenter}). To do this we will use  \cite[Corollary 5.2]{StepUniloc}, where it was proved that there exists $N\in \N$ such that for any commutative ring $R$, any $g\in G_{\P}(\Phi,R)$, and any $\beta\in \Phi$ the element $x_{\beta}(1)^g$ is a product of at most $N$ elementary root elements. Therefore for any $g\in G_\P (\Phi,R)$ we can describe its action on the elements of the set~$M$ by a formula and show that two actions coincide if and only if these are the same elements in $PG_\P (\Phi,R)$. 

Now we just need to prove that our parameters from $M$ are regularly definable. The idea here is that we take some parameters $\{ g_\alpha\mid \alpha \in \Phi\}$ and postulate that by above methods we define for any~$\alpha$ subgroups $Y_\alpha$ with ring structures (and they are all isomorphic to some ring~$S$), then with the help of this~$S$ we describe (by the scheme with matrices and equations, by a single formula) the group $PG_\P (\Phi,S)$ and postulate that this is  our initial group itself and the elements $g_\alpha$ correspond to $\overline{x_\alpha}(1)$. According to the description of isomorphisms of central quotients of Chevalley groups the rings $R$ and $S$ are isomorphic.  The proof yields, additionally,  that the groups $PG_\P (\Phi,-)$ are defined by one formula and that our parameters are defined up to automorphisms.

\medskip

Now we will show the detailed proof. 

 The first lemma essentially says that the root subgroups are definable.

{\lem \label{ElementaryRootSubgDefine} 
Let $\beta\in\Phi$, then there is a $(|\Phi|+1)$-ary formula $\IsElUnip_{\beta}(g,g_{\alpha}\mid \alpha\in\Phi)$ such that for any commutative ring $R$ and any $g\in PG_{\P}(\Phi,R)$ we have
$$
\left[PG_{\P}(\Phi,R)\models  \IsElUnip_{\beta}(g,\overline x_{\alpha}(1)\mid \alpha\in\Phi)\right]\Longleftrightarrow \left[\exists \xi\in R\mid g=\overline x_{\beta}(\xi)\right].
$$
}

\begin{proof}

Follows from the proofs of \cite[Propositions 3 and 4]{BMPDipphantine} in case of adjoint group. Note that the difference between $G_{\ad}(\Phi,R)$ and  $PG_{\P}(\Phi,R)$ does not affect these proofs.
\end{proof}

The next lemma allows to realize ring operations in terms of elementary unipotents.

{\lem \label{AddAndMult} Let $\beta\in\Phi$, then there are $(|\Phi|+3)$-ary formulas $\Add_{\beta}(g_1,g_2,g_3,g_{\alpha}\mid \alpha\in\Phi)$ and $\Mult_{\beta}(g_1,g_2,g_3,g_{\alpha}\mid \alpha\in\Phi)$ such that for any commutative ring $R$ and any $\xi_1$,$\xi_2$,$\xi_3\in R$ we have
$$
\left[PG_{\P}(\Phi,R)\models  \Add_{\beta}(\overline{x_{\beta}}(\xi_1),\overline{x_{\beta}}(\xi_2),\overline{x_{\beta}}(\xi_3),\overline{x_{\alpha}}(1)\mid \alpha\in\Phi)\right]\Longleftrightarrow \left[\xi_1+\xi_2=\xi_3\right]\tc
$$
and
$$
\left[PG_{\P}(\Phi,R)\models  \Mult_{\beta}(\overline{x_{\beta}}(\xi_1),\overline{x_{\beta}}(\xi_2),\overline{x_{\beta}}(\xi_3),\overline{x_{\alpha}}(1)\mid \alpha\in\Phi)\right]\Longleftrightarrow \left[\xi_1\xi_2=\xi_3\right]\tp
$$
}

\begin{proof}
    Follows from the proof of \cite[Theorem 4]{BMPDipphantine}.
\end{proof}

As we mentioned above, let us fix some faithful representation of $G_{\P}(\Phi,-)$, i.e a closed embedding of algebraic group schemes $G_{\P}(\Phi,-)\to \GL(n,-)$ for some $n$ (see \cite{ChevalleySemiSimp}); and let us fix some root $\alpha_0\in\Phi$.

{\lem \label{DecompConj}
Let $N$ be a positive integer, and let $\beta,\beta_1,\ldots,\beta_N\in\Phi$. Then there exists an $N+n^2+|\Phi|$-ary formula $\DecompConj_{\beta,\beta_1,\ldots,\beta_N}(g_1,\ldots, g_N,g_{i,j}\mid 1\le i,j\le n, g_{\alpha}\mid \alpha\in\Phi)$ such that for any commutative ring $R$, any $\zeta_1,\ldots,\zeta_N\in R$ and any $h\in G_{\P}(\Phi,R)$ with the matrix of $h$ in our representation being $(\xi_{i,j})$, we have
   \begin{multline*}
\!\!\!\!\!\!PG_{\P}(\Phi,R)\models \\
\!\!\!\!\models\DecompConj_{\beta,\beta_1,\ldots,\beta_N}(\overline{x_{\beta_1}}(\zeta_1),\ldots,\overline{x_{\beta_N}}(\zeta_N), \overline{x_{\alpha_0}}(\xi_{i,j})\mid 1\le i,j\le n, \overline{x_{\alpha}}(1)\mid \alpha\in\Phi)\Longleftrightarrow\\ \Longleftrightarrow x_{\beta}(1)^h=x_{\beta_1}(\zeta_1)\ldots x_{\beta_N}(\zeta_N)\; \mathrm{(in }\; G_{\P}(\Phi,R)\mathrm{)}\tp
\end{multline*}
}
\begin{proof}
     Note that the equality $x_{\beta}(1)^h=x_{\beta_1}(\zeta_1)\ldots x_{\beta_N}(\zeta_N)$ is equivalent to certain polynomial equality in $R$ on $\zeta_1,\ldots,\zeta_N$, and~$\xi_{i,j}$. So, in order to construct such a formula, we can turn elements $x_{\beta_i}(\zeta_i)$ into $x_{\alpha_0}(\pm\zeta_i)$ via conjugation with an element of the extended Weyl group, which is expressed through elements $x_{\alpha}(1)$; and then we use Lemma~\ref{AddAndMult} in order to write the required polynomial equality in the language of groups.
\end{proof}

{\lem \label{MatrixUpToCenter} There is an $(|\Phi|+n^2+1)$-ary formula $\MatUpToCent(g,g_{i,j}\mid 1\le i,j\le n, g_{\alpha}\mid \alpha\in\Phi)$ such that for any commutative ring $R$ and any $g$, $g_{i,j}\in PG_{\P}(\Phi,R)$, we have $PG_{\P}(\Phi,R)\models \MatUpToCent(g,g_{i,j}\mid 1\le i,j\le n, \overline{x_{\alpha}}(1)\mid \alpha\in\Phi)$ if and only if there are $\xi_{i,j}\in R$ such that the following holds:
\begin{enumerate}
    \item $g_{i,j}=\overline{x_{\alpha_0}}(\xi_{i,j})$ for all $1\le i,j\le n$;

    \item The matrix $(\xi_{i,j})$ is in the image of $G_{\P}(\Phi,R)$;

    \item The element $g$ of the group $PG_{\P}(\Phi,R)$ can be represented in the group $G_{\P}(\Phi,R)$ by the matrix $(\xi_{i,j})$.
\end{enumerate}
}

\begin{proof}
    We build this formula as the conjunction of the following three parts.

    Firstly, we use Lemma~\ref{ElementaryRootSubgDefine} and take the conjunction of $\IsElUnip_{\alpha_0}(g_{i,j}, g_{\alpha}\mid \alpha\in\Phi)$ for all $1\le i,j\le n$. So we ensure that $g_{i,j}=\overline{x_{\alpha_0}}(\xi_{i,j})$ for some $\xi_{i,j}\in R$.

    Secondly, we use Lemma~\ref{AddAndMult} to construct a formula stating that the matrix $(\xi_{i,j})$ satisfies the polynomial equations that determine the image of $G_{\P}(\Phi,R)$ in $\GL(n,R)$. We can do it because polynomials with integer coefficients can be obtained by the addition and multiplication.

    Now by~\cite[Corollary 5.2]{StepUniloc} there exists $N\in \N$ such that for any commutative ring $R$, any $g\in G_{\P}(\Phi,R)$, and any $\beta\in \Phi$ the element $x_{\beta}(1)^g$ is a product of at most $N$ elementary root elements.

    The third and the last part of our formula will be as follows:
    \begin{dmath*}
   \bigwedge_{\beta\in\Phi}\bigvee_{\beta_1,\ldots,\beta_N\in\Phi}\left[\exists h_1,\ldots,h_N\ \left(\bigwedge_{i=1}^N \IsElUnip_{\beta_i}(h_i,g_{\alpha}\mid \alpha\in\Phi)\right)\wedge (g_{\beta}^g=h_1\ldots h_N)\wedge\\ \wedge \DecompConj_{\beta,\beta_1,\ldots,\beta_N}(h_1,\ldots,h_N, g_{i,j}\mid 1\le i,j\le n, g_{\alpha}\mid \alpha\in\Phi) \right]\tc
    \end{dmath*}
   where the formula $\DecompConj_{\beta,\beta_1,\ldots,\beta_N}(h_1,\ldots,h_N, g_{i,j}\mid 1\le i,j\le n, g_{\alpha}\mid \alpha\in\Phi)$ is from Lemma~\ref{DecompConj}.

   Therefore, we achieved that the element $g$ and the element represented by the matrix $(\xi_{i,j})$ differs by multiplying by such an element $c$ that it commutes with all the $\overline{x_{\alpha}}(1)$. Since $PG_{\P}(\Phi,R)\le G_{\ad}(\Phi,R)$, it follows then from \cite[Theorem~3]{BMPDipphantine} that $c$ is trivial.
\end{proof}

Now we prove Theorem~\ref{RegularTrivCent}.
\begin{proof}
Note that item (4) directly follows from item (3), so we will prove only items (1)--(3).

    We construct the formula $\ph$ as a conjunction of formulas that essentially says the following:
    \begin{enumerate}
        \item The center of~$G$ is trivial.
        
        \item The formulas $\Add_{\alpha_0}$ and $\Mult_{\alpha_0}$ define binary operation on the set of such elements~$g$ that $\IsElUnip_{\alpha_0}(g,g_{\alpha}\mid \alpha\in\Phi)$ holds; and this set is a ring with respect to these operations.

        \item The ring $R$ from the item~2 is good.

        \item For any tuple $g_{i,j}$, $1\le i,j\le n$ such that $\IsElUnip_{\alpha_0}(g_{i,j},g_{\alpha}\mid \alpha\in\Phi)$ holds and that satisfy the equations that define $G_{\P}(\Phi,-)$ (i.e. satisfy the formula that is obtained from these equations by rewriting it in the language of groups using Lemma~\ref{AddAndMult}) there exists a unique $g$ such that $\MatUpToCent(g,g_{i,j}\mid 1\le i,j\le n, g_{\alpha}\mid \alpha\in\Phi)$ holds. In other words, the formula $\MatUpToCent$ defines a map from the set of matrices that satisfy the equations to the underlying set of our group.

        \item For any $g$ there exists a tuple $g_{i,j}$, $1\le i,j\le n$ such that $\MatUpToCent(g,g_{i,j}\mid 1\le i,j\le n, g_{\alpha}\mid \alpha\in\Phi)$ holds. In other words, the map from the previous item is surjective.

        \item For any $g$,$g'$,$g''$,$g_{i,j}$,$g_{i,j}'$,$g_{i,j}''$, if $\MatUpToCent(g,g_{i,j}\mid 1\le i,j\le n, g_{\alpha}\mid \alpha\in\Phi)$, $\MatUpToCent(g',g'_{i,j}\mid 1\le i,j\le n, g_{\alpha}\mid \alpha\in\Phi)$ and $\MatUpToCent(g'',g''_{i,j}\mid 1\le i,j\le n, g_{\alpha}\mid \alpha\in\Phi)$ hold and the matrix $(g_{i,j})$ is the product of matrices $(g'_{i,j})$ and $(g''_{i,j})$ in the sense of operations defined by $\Add_{\alpha_0}$ and $\Mult_{\alpha_0}$, then $g=g'g''$. In other words, the map from the previous item is a group homomorphism.

        \item For any tuple  $g_{i,j}$, $1\le i,j\le n$ we have $\MatUpToCent(e,g_{i,j}\mid 1\le i,j\le n, g_{\alpha}\mid \alpha\in\Phi)$ holds if and only if it represents the matrix from the center of $G_{\P}(\Phi,R)$ (note that the center is defined by equations). In other words, the kernel of the homomorphism from the previous items is the center.

        \item The homomorphism from the previous items sends the elements $x_{\alpha}(1)$ to the elements $g_{\alpha}$ (the matrices $x_{\alpha}(1)$ are with integer entries, and any particular integer is definable in the language of rings).
    \end{enumerate}

    \smallskip

Now let us prove that the formula constructed above satisfies items (1)--(3) of the Theorem~\ref{RegularTrivCent}.

Firstly, by construction, for any good commutative ring $R$  we have 
$$
PG_{\P}(\Phi,R)\models \ph(\overline{x_{\alpha}}(1)\mid \alpha\in\Phi)\tp
$$
Therefore, the implication from the right to the left in item (1) of Theorem~\ref{RegularTrivCent} holds true because the automorphisms respects any logical formula.

Secondly, item~(3) of Theorem~\ref{RegularTrivCent} also holds by construction.

Now let us prove the implication from the left to the right in item (1). Assume that for some $g_{\alpha}\in PG_{\P}(\Phi,R)$ we have
$$
PG_{\P}(\Phi,R)\models \ph(g_{\alpha}\mid \alpha\in\Phi)\tp
$$
By item~(3), there exists an isomorphism $f:PG_\P (\Phi, S) \to PG_\P (\Phi,R)$ for some ring~$S$ such that $f(\overline{x_\alpha}(1))=g_\alpha$ for all $\alpha \in \Phi$. By Theorem~\ref{isomorphic} the (good) rings $R$ and $S$ are isomorphic. Now composing $f$ with the ring isomorphism from $PG_\P (\Phi,S)$ to $PG_\P (\Phi,R)$, we obtained the required automorphism.

It remains to prove the item~(2). By item~(1) it is enough to prove that there exists a bi-interpretation of $R$ and $PG_{\P}(\Phi,R)$ such that the interpretation of $PG_{\P}(\Phi,R)$ in $R$ is absolute and the parameters for interpretation of $R$ in $PG_{\P}(\Phi,R)$ are $\overline{x_{\alpha}}(1)$.

We interpret $PG_{\P}(\Phi,R)$ in $R$ by an obvious way: elements of $PG_{\P}(\Phi,R)$ are matrices over $R$ that satisfy certain equations; equiality of elements is an equivalence up to the center; and the multiplication is the multiplication of matrices.

Further we interpret $R$ in $PG_{\P}(\Phi,R)$ using Lemmas~\ref{ElementaryRootSubgDefine} and~\ref{AddAndMult}.

It is easy to see that the set of pairs $(\xi,g)\in R\times G_{\P}(\Phi,R)$ such that $g$ differs from $x_{\alpha_0}(\xi)$ by a central element is a set of points of a closed subscheme in $\mathbb{A}_1\times G_{\P}(\Phi,-)$; hence the isomorphism between $R$ and its interpretation inside $PG_{\P}(\Phi,R)$ is defined by polynomial formulas.

The isomorphism between $PG_{\P}(\Phi,R)$ and its interpretation inside $R$ can be defined using Lemma~\ref{MatrixUpToCenter}.
\end{proof}

{\cor \label{Corollary_Adjoint}
Let $\Phi$ be an irreducible root system of rank $>1$, $R$ be a good commutative ring. Then the adjoint Chevalley group $G_{ad}(\Phi,R)$ satisfies all conclusions of Theorem~\ref{RegularTrivCent}.

In particular, if an arbitrary group $G$ is elementary equivalent to $G_{ad}(\Phi,R)$, then there exists a good ring $S$ such that $G\cong G_{ad}(\Phi,S)$.
}

\begin{proof}
 This corollary follows from Theorem~\ref{RegularTrivCent} and the fact that for any adjoint Chevalley group its center is trivial (see~\cite{AbeHurley}).  
\end{proof}

{\cor \label{Finite_Axiomatize}
Let $\Phi$ be an irreducible root system of rank $>1$, $\P$ be any corresponding weight lattice, $\mathbf G_{\Phi,\P}$ and $\mathbf{PG}_{\Phi,\P}$ be the classes of 
\begin{enumerate}
\item all Chevalley groups $G_\P (\Phi, R)$ with trivial centers and groups $PG_\P (\Phi, R)$ over arbitrary commutative rings for $\Phi = A_\ell, D_\ell, E_\ell$, $\ell \geqslant 3$;
\item all Chevalley groups $G_\P (\Phi, R)$ with trivial centers and groups $PG_\P (\Phi, R)$ over arbitrary commutative rings with $1/2$ for $\Phi = A_2, B_\ell, C_\ell, F_4$, $\ell \geqslant 2$;   
\item all Chevalley groups $G_\P (\Phi, R)$  over arbitrary commutative rings with $1/2$ and $1/3$ for $\Phi = G_2$. 
\end{enumerate}
Then the classes $\mathbf G_{\Phi,\P}$  and $\mathbf{PG}_{\Phi,\P}$ are finitely axiomatizable.
}

\begin{proof}
Directly follows from Theorem~\ref{RegularTrivCent}.    
\end{proof}

{\cor\label{Axiomatize} Let $\Phi$ be an irreducible root system of rank $>1$, $\P$ be any corresponding weight lattice. Let $T$ be a collection (possibly infinite) of first order sentences in the language of rings such that $T$ includes the sentence "$1/2\in R$" in case $\Phi = A_2, B_\ell, C_\ell, F_4$, $G_2$ $\ell \geqslant 2$, and $T$ includes the sentence "$1/3\in R$" in case $\Phi = G_2$. Then the class of groups that are isomorphic to $PG_\P (\Phi, R)$ for some ring $R$ that satisfy all sentences from $T$ is axiomatizable, i.e. it can be described by (possibly infinite) collection of first order sentences.
}

\begin{proof}
   In case where $T$ consist of one sentence, the statement follows directly from Theorem~\ref{RegularTrivCent}; moreover in this case the class in question can be describe by one sentence. In general case we take each sentence in $T$, construct the corresponding sentence for groups and take the collection of all such sentences. Clearly, the groups in question must satisfy all these sentences. Conversely, if a group $G$ satisfies all these sentences, then for every sentence from $T$ there exists a ring $R$ such that $G\simeq PG_\P (\Phi, R)$ and $R$ satisfy this sentence; but by Theorem~\ref{isomorphic} all these rings must be isomorphic; hence, they satisfy all the sentences from $T$. 
\end{proof}

\begin{exmp}
    The class of groups that are isomorphic to $PG_\P (\Phi, K)$, where $K$ is a pseudo-finite field (of admissible characteristic) is axiomatizable (for details about logical properties of groups over pseudo-finite fields see \cite{HrushPill}).
\end{exmp}

{\cor\label{QFA} Let $\Phi$ be an irreducible root system of rank $>1$, $R$ be a finitely generated commutative ring. Assume that the group $PG_\P (\Phi, R)$ is finitely generated. Then this group is QFA. In particular, all the finitelly generated adjoint Chevalley groups of rank $>1$ are QFA, cf.~\cite{Segal-Tent}. In particular these groups are first order rigid.}
\begin{proof}
 By~\cite{Aschenbrenner} all finitely generated commutative rings are QFA. Let $\eta$ be a proposition that defines $R$ among all finitely generating rings. Now construct the formula  $\ph(g_{\alpha}\mid \alpha\in\Phi)$ in the same way that we did in the proof of Theorem~\ref{RegularTrivCent} except we do not require the ring obtained from the group to be good and add the item that says that the ring satisfies $\eta$.

We claim that the proposition $\exists \{g_{\alpha}\}_{\alpha\in\Phi}\; \ph(g_{\alpha}\mid \alpha\in\Phi)$ defines $PG_\P (\Phi, R)$ among all finitely generating groups.

Indeed, let $G$ be a finitely generated group that satisfies that proposition. By construction of the formula $\ph$ we have $G\simeq PG_\P (\Phi,S)$, where $S$ satisfies $\eta$ (note that the proof of this part of Theorem~\ref{RegularTrivCent} does not require the ring to be good). Since $G$ is a finitely generated group it follows that $S$ is a finitely generated ring. Hence we have $R\simeq S$. Therefore, we have $G\simeq PG_\P (\Phi,R)$.
\end{proof}

\section{Proof of Theorem~\ref{RegularBoundGen}}

{\lem \label{ElementaryRootSubgDefine2} 
Let $\beta\in\Phi$, then there is a $(|\Phi|+1)$-ary formula $\IsElUnip_{\beta}(g,g_{\alpha}\mid \alpha\in\Phi)$ such that for any commutative ring $R$ such that $1/2\in R$ in case $\Phi=C_2$ and any $g\in G_{\P}(\Phi,R)$ we have
$$
\left[G_{\P}(\Phi,R)\models  \IsElUnip_{\beta}(g,x_{\alpha}(1)\mid \alpha\in\Phi)\right]\Leftrightarrow \left[\exists \xi\in R\mid g=x_{\beta}(\xi)\right].
$$
}

\begin{proof}
     Follows from the proofs of \cite[Propositions 3 and 4]{BMPDipphantine}.
\end{proof}

{\lem \label{AddAndMult2} Let $\beta\in\Phi$, then there is a $(|\Phi|+3)$-ary formulas $\Add_{\beta}(g_1,g_2,g_3,g_{\alpha}\mid \alpha\in\Phi)$ and $\Mult_{\beta}(g_1,g_2,g_3,g_{\alpha}\mid \alpha\in\Phi)$ such that for any commutative ring $R$ such that $1/2\in R$ in case $\Phi=C_2$ and any $\xi_1$,$\xi_2$,$\xi_3\in R$ we have
$$
\left[G_{\P}(\Phi,R)\models  \Add_{\beta}(x_{\beta}(\xi_1),x_{\beta}(\xi_2),x_{\beta}(\xi_3),x_{\alpha}(1)\mid \alpha\in\Phi)\right]\Leftrightarrow \left[\xi_1+\xi_2=\xi_3\right]\tc
$$
and
$$
\left[G_{\P}(\Phi,R)\models  \Mult_{\beta}(x_{\beta}(\xi_1),x_{\beta}(\xi_2),x_{\beta}(\xi_3),x_{\alpha}(1)\mid \alpha\in\Phi)\right]\Leftrightarrow \left[\xi_1\xi_2=\xi_3\right]\tp
$$
}

\begin{proof}
    Follows from the proof of \cite[Theorem 4]{BMPDipphantine}.
\end{proof}

{\lem \label{MatrixUpToCenter2} There is an $(|\Phi|+n^2+1)$-ary formula $\MatUpToCent(g,g_{i,j}\mid 1\le i,j\le n, g_{\alpha}\mid \alpha\in\Phi)$ such that for any commutative ring $R$ such that $1/2\in R$ in case $\Phi=C_2$ and any $g$, $g_{i,j}\in G_{\P}(\Phi,R)$, we have $G_{\P}(\Phi,R)\models \MatUpToCent(g,g_{i,j}\mid 1\le i,j\le n, x_{\alpha}(1)\mid \alpha\in\Phi)$ if and only if there are $\xi_{i,j}\in R$ such that the following holds:
\begin{enumerate}
    \item $g_{i,j}=x_{\alpha_0}(\xi_{i,j})$ for all $1\le i,j\le n$;

    \item The matrix $(\xi_{i,j})$ is in the image of $G_{\P}(\Phi,R)$ under the representation $\vpi$;

    \item The element $g$ differs from the element represented by the matrix $(\xi_{i,j})$ by multiplication by an element from the center of $G_{P}(\Phi,R)$.
\end{enumerate}
}

\begin{proof}
    Similar to the proof of Lemma~\ref{MatrixUpToCenter}.
\end{proof}

{\lem \label{Matrix} There is an $|\Phi|+(l+1)(n^2+1)$-ary formula $\Matrix(g, g_{i,j}\mid 1\le i,j\le n, g_{\alpha}\mid \alpha\in\Phi;\; g_k\mid 1\le k\le l;\; g_{ij}^{(k)}\mid 1\le k\le l,\; 1\le i,j\le n)$ such that for any commutative ring $R$ such that $1/2\in R$ in case $\Phi=C_2$ and any $g$, $g_{i,j}$,$g_k\in G_{\P}(\Phi,R)$ such that the group $G_{\P}(\Phi,R))$ is $N$-boundedly generated by the set $X_{g_1,\ldots,g_l}$, we have $G_{\P}(\Phi,R)\models \Matrix(g, g_{i,j}\mid 1\le i,j\le n, x_{\alpha}(1)\mid \alpha\in\Phi;\; g_k\mid 1\le k\le l;\; x_{\alpha_0}(\vpi(g_k)_{ij})\mid 1\le k\le l,\; 1\le i,j\le n)$ if and only if there are $\xi_{i,j}\in R$ such that the following holds:

\begin{enumerate}
    \item $g_{i,j}=x_{\alpha_0}(\xi_{i,j})$ for all $1\le i,j\le n$;

    \item The matrix $(\xi_{i,j})$ equals $\vpi(g)$.
\end{enumerate}}

\begin{proof}
 We build this formula as the conjunction of two parts.

    In the first part, we use Lemma~\ref{ElementaryRootSubgDefine2} and take the conjunction of $\IsElUnip_{\alpha_0}(g_{i,j}, g_{\alpha}\mid \alpha\in\Phi)$ for all $1\le i,j\le n$. So we ensure that $g_{i,j}=x_{\alpha_0}(\xi_{i,j})$ for some $\xi_{i,j}\in R$.

    Now let $\Omega$ be the set of all words of length $N$ in letters $\{c,d,x_1,\ldots,x_l\}$. To each $\omega\in\Omega$ we assign a term in the language of groups build as follows: we take two variables for each letter $c$ and one variable for each letter $d$, all variables are distinct from each other and from the letters $x_1$,$\ldots$,$x_l$. Now we replace each letter $c$ with the commutator of corresponding variables, and each letter $d$ with the $m$-th power of the corresponding variable. 

    The second part is constructed as a disjunction by the set $\Omega$ of the formulas that says the following:

\begin{enumerate}
    \item the element $g$ can be represented by the term constructed from the corresponding element $\omega\in\Omega$ as described above with elements $g_1$,$\ldots$,$g_l$ being substituted as $x_1$,$\ldots$,$x_l$;

    \item the elements $g_{ij}$ satisfy the formula constructed using Lemmas~\ref{ElementaryRootSubgDefine2} and~\ref{AddAndMult2} that essentially says that the matrix $(\xi_{ij})$ can be represented by the same term as above  with matrices $(\xi_{ij}^{(1)})$,$\ldots$,$(\xi_{ij}^{(l)})$ being substituted as $x_1$,$\ldots$,$x_l$, assuming $g_{ij}^{(k)}=x_{\alpha_0}(\xi_{ij}^{(k)})$ for $1\le k\le l$;

    \item for any variable that comes from the letter $c$ or $d$ in the corresponding element $\omega\in\Omega$ the corresponding element used in item~1 differs from the element represented by the corresponding matrix used in item~2 by an element from the center of $G_{P}(\Phi,R)$ (this item uses Lemma~\ref{MatrixUpToCenter2}).
\end{enumerate}

Note that if $a_1$ differs from $a_2$ by a central element and $b_1$ differs from $b_2$ by a central element, then $[a_1,b_1]=[a_2,b_2]$ and $a_1^m=a_2^m$. Therefore, our formula has the required property.
\end{proof}

Now we prove Theorem~\ref{RegularBoundGen}.

\begin{proof}
Note that item (4) directly follows from items (1)--(3), so we will prove only items (1)--(3).

    We construct the desired formula $\ph$ as a conjunction of formulas that essentially says the following:
    \begin{enumerate}
\item $G$ is $N$-boundedly generated by the set $X^{(m)}_{g_1,\ldots,g_l}$.

        \item The formulas $\Add_{\alpha_0}$ and $\Mult_{\alpha_0}$ define binary operation on the set of such elements~$g$ that $\IsElUnip_{\alpha_0}(g;g_{\alpha}\mid \alpha\in\Phi)$ holds; and this set is a ring with respect to these operations.

        \item The ring $R$ from the item~2 is good.
        
         \item For any tuple $g_{i,j}$, $1\le i,j\le n$ such that $\IsElUnip_{\alpha_0}(g_{i,j};g_{\alpha}\mid \alpha\in\Phi)$ holds and that satisfy the equations that define $G_{\P}(\Phi,-)$ (i.e. satisfy the formula that is obtained from these equations by rewriting it in the language of groups using Lemma~\ref{AddAndMult2}) there exists a unique $g$ such that $\Matrix(g; g_{i,j}\mid 1\le i,j\le n, g_{\alpha}\mid \alpha\in\Phi;\; g_k\mid 1\le k\le l;\; g_{ij}^{(k)}\mid 1\le k\le l,\; 1\le i,j\le n)$ holds.

        \item For any $g$ there exists a unique tuple $g_{i,j}$, $1\le i,j\le n$ such that it satisfy the equations as in item~(4) and $\Matrix(g; g_{i,j}\mid 1\le i,j\le n, g_{\alpha}\mid \alpha\in\Phi;\; g_k\mid 1\le k\le l;\; g_{ij}^{(k)}\mid 1\le k\le l,\; 1\le i,j\le n)$ holds.

        \item For any $g$,$g'$,$g''$,$g_{i,j}$,$g_{i,j}'$,$g_{i,j}''$, if $\Matrix(g; g_{i,j}\mid 1\le i,j\le n, g_{\alpha}\mid \alpha\in\Phi;\; g_k\mid 1\le k\le l;\; g_{ij}^{(k)}\mid 1\le k\le l,\; 1\le i,j\le n)$, $\Matrix(g'; g_{i,j}'\mid 1\le i,j\le n, g_{\alpha}\mid \alpha\in\Phi;\; g_k\mid 1\le k\le l;\; g_{ij}^{(k)}\mid 1\le k\le l,\; 1\le i,j\le n)$ and $\Matrix(g''; g_{i,j}''\mid 1\le i,j\le n, g_{\alpha}\mid \alpha\in\Phi;\; g_k\mid 1\le k\le l;\; g_{ij}^{(k)}\mid 1\le k\le l,\; 1\le i,j\le n)$ hold and the matrix $(g_{i,j})$ is the product of matrices $(g'_{i,j})$ and $(g''_{i,j})$ in the sense of operations defined by $\Add_{\alpha_0}$ and $\Mult_{\alpha_0}$, then $g=g'g''$.

        \item In the correspondence between elements and matrices defined by items~3 and~4, matrices of elements $x_{\alpha}(1)$ correspond to the elements $g_{\alpha}$ (these are matrices with integer entries, and any particular integer is definable in the language of rings); and for any $1\le k\le l$ the matrix represented by the tuple $(g_{i,j}^{(k)})$ correspond to the element~$g_k$.
    \end{enumerate}

    Now let us prove that the formula constructed above satisfies items (1)--(3) of the Theorem~\ref{RegularBoundGen}. 

Firstly, the implication from the right to the left in item (1) holds true by construction of $\ph$.

Secondly, item~(3) also holds by construction.

Now let us prove the implication from the left to the right in item (1). Assume that for some $g_{\alpha},g_k,g_{ij}^{(k)}\in G_{\P}(\Phi,R)$ we have
$$
G_{\P}(\Phi,R)\models \ph(g_{\alpha}\mid \alpha\in\Phi;\; g_k\mid 1\le k\le l;\; g_{ij}^{(k)}\mid 1\le k\le l,\; 1\le i,j\le n)\tp
$$
By item~(3), there exists an isomorphism $f\colon G_\P (\Phi, S) \to G_\P (\Phi,R)$ for some ring~$S$ such that

$\bullet$ $f(g_{\alpha})= x_{\alpha}(1)$ for all $\alpha\in\Phi$; 

$\bullet$ $f(g_{ij}^{(k)})=x_{\alpha_0}(\xi_{ij}^{(k)})$ for all $1\le k\le l$ for some $\xi_{ij}^{(k)}\in R$;

$\bullet$ $\vpi(f(g_k))$ is the matrix $(\xi_{ij}^{(k)})$ for all $1\le k\le l$.

By Theorem~\ref{isomorphic} the rings $R$ and $S$ are isomorphic (here we use that they are good). Now composing $f$ with the ring isomorphism from $G_\P (\Phi,S)$ to $G_\P (\Phi,R)$, we obtain the required automorphism.

It remains to prove the item~(2). By item~(1) it is enough to prove that there exists a bi-interpretation of $R$ and $G_{\P}(\Phi,R)$ such that the interpretation of $G_{\P}(\Phi,R)$ in $R$ is absolute and the parameters for interpretation of $R$ in $G_{\P}(\Phi,R)$ are $x_{\alpha}(1)$,$g_k$, and $x_{\alpha_0}(\xi_{ij}^{(k)})$, where $\vpi(g_k)=(\xi_{ij}^{(k)})$.

We interpret $G_{\P}(\Phi,R)$ in $R$ by an obvious way: elements of $G_{\P}(\Phi,R)$ are matrices over $R$ that satisfy certain equations and the multiplication is the multiplication of matrices.

Further we interpret $R$ in $G_{\P}(\Phi,R)$ using Lemmas~\ref{ElementaryRootSubgDefine2} and~\ref{AddAndMult2}.

It is easy to see that the set of pairs $(\xi,g)\in R\times G_{\P}(\Phi,R)$ such that $\vpi(g)=x_{\alpha_0}(\xi)$ is a set of points of a closed subscheme in $\mathbb{A}_1\times G_{\P}(\Phi,-)$; hence the isomorphism between $R$ and its interpretation inside $G_{\P}(\Phi,R)$ is defined by polynomial formulas.

The isomorphism between $G_{\P}(\Phi,R)$ and its interpretation inside $R$ can be defined using Lemma~\ref{Matrix}.
\end{proof}

{\cor\label{Axiomatize2} Let $\Phi$ be an irreducible root system of rank $>1$, $\P$ be any corresponding weight lattice, and $l$,$N$ be positive integers. Let $T$ be a collection (possibly infinite) of first order sentences in the language of rings such that $T$ includes the sentence "$1/2\in R$" in case $\Phi = A_2, B_\ell, C_\ell, F_4$, $G_2$ $\ell \geqslant 2$, and $T$ includes the sentence "$1/3\in R$" in case $\Phi = G_2$, and such that $T$ includes the sentence "there is a $l$-tuple $g_1,\ldots,g_l\in G$ such that the group $G$ is $N$-boundednly generated by the set $X_{g_1,\ldots,g_l}$" . Then the class of groups that are isomorphic to $G_\P (\Phi, R)$ for some ring $R$ that satisfy all sentences from $T$ is axiomatizable, i.e. can be described by (possibly infinite) collection of first order sentences.
}

\begin{proof}
   Similar to Corollary~\ref{Axiomatize}.
\end{proof}

\begin{exmp}
    The class of groups that are isomorphic to $G_\P (\Phi, K)$, where $K$ is a pseudo-finite field (of admissible characteristic) is axiomatizable.
\end{exmp}

{\cor\label{QFA2} Let $\Phi$ be an irreducible root system of rank $>1$, $R$ be a finitely generated commutative ring. In the case $\Phi=C_2$ assume that $1/2\in R$. Assume that the group $G_\P (\Phi, R)$ is finitely generated, and boundedly generated by some set $X^{(m)}_{g_1,\ldots,g_l}$ (for example that is true, when $R$ is a Dedekind domain of arithmetic type). Then this group is QFA. In particular these groups are first order rigid.}

\begin{proof}
   Similar to Corollary~\ref{QFA}.
\end{proof}

{\rem We believe that the assumption $1/2\in R$ for $\Phi=C_2$ can be lifted, but this case requires additional arguments.}

\section{A counterexample for the infinite center}

In this section we will give the example of a Chevalley group with infinite center, which is not bi-interpretable with the corresponding ring and which is elementary equivalent to a non-Chevalley group.

\medskip

Consider the Chevalley group $G_{\P}(A_3,R)$, where $\P$ is the only suitable lattice that is neither root lattice nor weight lattice (note that such lattices are possible only for systems of types $A$ and $D$). That is we take the special orthogonal group $\SO(6,R)$, i.e. the group of $6\times 6$ matrices that preserves the quadratic form $x_1x_6+x_2x_5+x_3x_4$, and have the trivial Dickson invariant. The corresponding elementary subgroup will be denoted by $EO(6,R)$.

The torus $TO(6,R)\le \SO(6,R)$ consists of diagonal matrices that have shape $\diag(\xi_1,\xi_2,\xi_3, \xi_3^{-1},\xi_2^{-1},\xi_1^{-1})$. where $\xi_i\in R^*$.

{\lem\label{TcapE} The element $\diag(\xi_1,\xi_2,\xi_3, \xi_3^{-1},\xi_2^{-1},\xi_1^{-1})$ belong to $EO(6,R)$ if and only if $\xi_1\xi_2\xi_3\in (R^*)^2$.}
\begin{proof}
    The simply connected Chevalley group $G_{\mathrm{sc}}(A_3,R)$ is the group $\SL(4,R)$. The element from $TO(6,R)$ belong to $EO(6,R)$ if and only if it belong to the image of the toric subgroup of $\SL(4,R)$ under the natural map ${\SL(4,R)}\to{\SO(6,R)}$. Indeed, the toric elements of $\SL(4,R)$ are elementary; hence so are their images. Conversely, every elementary element belongs to the image of $\SL(4,R)$ and the preimage of $TO(6,R)$ is the torus.

    Now the map on tori is given by
    $$
    \diag(\zeta_1,\zeta_2,\zeta_3,\zeta_4)\mapsto \diag(\zeta_1\zeta_2,\zeta_1\zeta_3,\zeta_1\zeta_4,\zeta_3\zeta_4,\zeta_2\zeta_4,\zeta_2\zeta_3)\tp
    $$

    Since $\zeta_1\zeta_2\zeta_3\zeta_4=1$, it follows that the product of the first three diagonal entries of the element from the image is equal to $\zeta_1^2\in (R^*)^2$. Conversely, suppose we have $\xi_1\xi_2\xi_3=\zeta_1^2$ for some $\zeta_1\in R^*$, Then we set $\zeta_2=\xi_1\zeta_1^{-1}$, $\zeta_3=\xi_2\zeta_1^{-1}$, $\zeta_4=\xi_3\zeta_1^{-1}$; hence we have $\zeta_1\zeta_2\zeta_3\zeta_4=1$. and 
    $$
    \diag(\xi_1,\xi_2,\xi_3,\xi_3^{-1},\xi_2^{-1},\xi_1^{-1})=\diag(\zeta_1\zeta_2,\zeta_1\zeta_3,\zeta_1\zeta_4,\zeta_3\zeta_4,\zeta_2\zeta_4,\zeta_2\zeta_3)\tp
    $$
\end{proof}

{\lem\label{ForFt} We have
$$
\SO(6,\F_2[t])=EO(6,\F_2[t])
$$}

\begin{proof}
    Follows from~\cite{StavSerreConj}.
\end{proof}

{\prop \label{Product} Set $R=\F_2[t,\eps]/(\eps^2)$. Then
$$
\SO(6,R)=EO(6,R)\times (\Z/2\Z)^{\aleph_0}\tp
$$}
\begin{proof}
The center of $\SO(6,R)$ consists of scalar matrices with the scalar being the square root of one. Thus the center is isomorphic to $\mu_2(R)=(\Z/2\Z)^{\aleph_0}$. Clearly, the center commutes with $EO(6,R)$. Therefore, it remains to show that $EO(6,R)\cap Z(\SO(6,R))=e$, and that $\SO(6,R)$ is generated by $EO(6,R)$ and $Z(\SO(6,R))$.

Let us prove the first claim. Let $g\in EO(6,R)\cap Z(\SO(6,R))$. Since $g\in Z(\SO(6,R))$, we have $g=\diag(\xi,\ldots,\xi)$, with $\xi^2=1$. By Lemma~\ref{TcapE} we have $\xi^3\in (R^*)^2=\{1\}$; hence $\xi=\xi^3=1$.

Let us prove the second claim. Clearly, $\SO(6,R)$ is generated by $\SO(6,\F_2[t])$ and the congruence subgroup $\SO(6,\eps R,R)$. By Lemma~\ref{ForFt} the group $\SO(6,\F_2[t])$ is contained in the elementary subgroup. Since $\eps$ is nilpotent, it follows that the congruence subgroup $\SO(6,\eps R,R)$ is contained in the big Gauss cell; hence it is contained in $TO(6,R)EO(6,R)$. Therefore, it remains to show that any element of $TO(6,R)$ is a product of an elementary element and a central element. Note that all the elements of $R^*$ are square roots of one. So given an element $\diag(\xi_1,\xi_2,\xi_3, \xi_3,\xi_2,\xi_1)\in TO(6,R)$ we can multiply it by the central element $\diag(\xi_1\xi_2\xi_3,\ldots \xi_1\xi_2\xi_3)$, and by Lemma~\ref{TcapE} the result will be in the elementary subgroup.
\end{proof}

{\cor\label{NotBiinterpretable} For $R=\F_2[t,\eps]/(\eps^2)$ the group $\SO(6,R)$ is not bi-interpretable with $R$.}
\begin{proof}
    Assume the converse. Then it is easy to see that the subgroup of $\Aut(\SO(6,R))$ that stabilises the corresponding parameters must be isomorphic to a certain subgroup of $\Aut(R)$. However, the group $\Aut(R)$ is countable (because $R$ is countable and the automorphism is uniquely defined by the images of $t$ and $\eps$); and it follows from Proposition~\ref{Product} that for any finite subset of $\SO(6,R)$ the subgroup of $\Aut(\SO(6,R))$ that stabilises the elements of this set is uncountable. This is a contradiction.
\end{proof}

{\cor\label{NonStandartModel} For $R=\F_2[t,\eps]/(\eps^2)$ there is the group $G$ that is elementary equivalent to $\SO(6,R)$, but not isomorphic to any Chevalley group over any ring.}
\begin{proof}
    Set $G=EO(6,R)\times (\Z/2\Z)^I$, where $I=[0,1]$. It follows from Proposition~\ref{Product} that $G$ is elementary equivalent to $\SO(6,R)$.

    It is easy to see that an infinite Chevalley group over a ring must have the same cardinality as its derived subgroup. Therefore, $G$ is not a Chevalley group.
\end{proof}

\medskip

{\bf Acknowledgements.}
Our sincere thanks go to Eugene Plotkin for very useful discussions regarding various aspects of this work and permanent attention to it. 

\bigskip


\end{document}